\documentclass{article}
\usepackage{amssymb}
\usepackage{latexsym}

\usepackage{graphicx}

\title{On hypergraph Lagrangians}

\author{ Qingsong Tang \thanks{College of Sciences, Northeastern University, Shenyang, 110819, China and School of Mathematics, Jilin University, Changchun 130012, P.R. China. Email: t\_qsong@sina.com.cn}  \and   Xiaojun Lu \thanks{College of Sciences, Northeastern University, Shenyang, 110819, China. Email: luxiaojun0625@sina.com}   \and  Xiangde Zhang \thanks{College of Sciences, Northeastern University, Shenyang, 110819, China. Email: zhangxdneu@163.com} \and Cheng Zhao \thanks{Department of Mathematics and Computer Science, Indiana State University, Terre Haute, IN, 47809 and School of Mathematics, Jilin University, Changchun 130012, P.R. China. Email: cheng.zhao@indstate.edu} }

\date{}
\newtheorem{defi}{Definition}[section]
\newtheorem{theo}{Theorem}[section]

\newtheorem{remark}[theo]{Remark}

\newtheorem{lemma}[theo]{Lemma}

\newtheorem{coro}[theo]{Corollary}
\newtheorem{con}[theo]{Conjecture}

\newtheorem{fact}[theo]{Fact}

\newcommand{\qed}{\hspace*{\fill} \rule{7pt}{7pt}}

\begin{document}
\maketitle
\begin{abstract}
It is conjectured by Frankl and F\"{u}redi that the $r$-uniform hypergraph with $m$ edges formed by taking the first $m$ sets in the colex ordering of ${\mathbb N}^{(r)}$ has the largest Lagrangian of all $r$-uniform hypergraphs with $m$ edges in \cite{FF}. Motzkin and Straus' theorem confirms this conjecture when $r=2$. For $r=3$, it is shown by Talbot  in \cite{T} that this conjecture is true when $m$ is in certain ranges.  In this paper,  we explore the connection between the clique number and Lagrangians for $r$-uniform hypergraphs. As an implication of this connection,  we prove that  the $r$-uniform hypergraph with $m$ edges formed by taking the first $m$ sets in the colex ordering of ${\mathbb N}^{(r)}$ has the largest Lagrangian of all $r$-uniform graphs with $t$ vertices and $m$ edges satisfying ${t-1\choose r}\leq m \leq {t-1\choose r}+ {t-2\choose r-1}-[(2r-6)\times2^{r-1}+2^{r-3}+(r-4)(2r-7)-1]({t-2\choose r-2}-1)$ for $r\geq 4.$
\end{abstract}

\\

\section{Introduction}

For a set $V$ and a positive integer $r$, let  $V^{(r)}$ be  the family of all $r$-subsets of $V$. An $r$-uniform hypergraph or $r$-graph $G$ consists of a set $V(G)$ of vertices and a set $E(G) \subseteq V(G) ^{(r)}$ of edges. When $r=2$, an $r$-graph is a simple graph.  When $r\ge 3$,  an $r$-graph is often called a hypergraph. An edge $e=\{a_1, a_2, \ldots, a_r\}$ will be simply denoted by $a_1a_2 \ldots a_r$. Let $K^{(r)}_t$ denote the complete $r$-graph on $t$ vertices, that is the $r$-graph on $t$ vertices containing all possible edges. A complete $r$-graph on $t$ vertices is also called a clique with order $t$.  A clique is said to be maximum if it has maximum cardinality.   Let ${\mathbb N}$ be the set of all positive integers.  For an integer $n \in {\mathbb N}$, let $[n]$ denote the set $\{1, 2, 3, \ldots, n\}$. Let $[n]^{(r)}$  represent the  complete $r$-graph on the vertex set $[n]$.

For an $r$-graph $G=(V,E)$, denote the $(r-1)$-neighborhood of a vertex $i \in V$ by $E_i=\{A \in V^{(r-1)}: A \cup \{i\} \in E\}$. Similarly,  denote the $(r-2)$-neighborhood of a pair of vertices $i,j \in V$ by $E_{ij}=\{B \in V^{(r-2)}: B \cup \{i,j\} \in E\}$.  Denote the complement of $E_i$ by $E^c_i=\{A \in V^{(r-1)}: A \cup \{i\} \in V^{(r)} \backslash E\}$. Also, denote the complement of $E_{ij}$ by
$E^c_{ij}=\{B \in V^{(r-2)}: B \cup \{i,j\} \in V^{(r)} \backslash E\}$ and $E_{i\setminus j}=E_i\cap E^c_j.$

\begin{defi}
For  an $r$-graph $G=([n], E(G))$ and a vector $\vec{x}=(x_1,\ldots,x_n) \in R^n$,
define
$$\lambda (G,\vec{x})=\sum_{i_1i_2 \cdots i_r \in E(G)}x_{i_1}x_{i_2}\ldots x_{i_r}.$$
Let $S=\{\vec{x}=(x_1,x_2,\ldots ,x_n): \sum_{i=1}^{n} x_i =1, x_i
\ge 0 {\rm \ for \ } i=1,2,\ldots , n \}$. The Lagrangian\footnote{Let us note that this use of the name Lagrangian is at odds with the tradition. Indeed, names as Laplacian, Hessian, Gramian, Grassmanian,  etc., usually denote a structured object like matrix,  operator, or manifold, and not just a single  number.}  of
$G$, denoted by $\lambda (G)$, is the maximum
 of the above homogeneous  function over the standard simplex $S$. Precisely,
 $$\lambda (G) = \max \{\lambda (G, \vec{x}): \vec{x} \in S \}.$$
\end{defi}

The value $x_i$ is called the  weight  of the vertex $i$.
A vector $\vec{x}=(x_1, x_2, \ldots, x_n) \in R^n$ is called a feasible weighting for $G$ if
$\vec{x}\in S$. A vector $\vec{y}\in S$ is called an   optimal weighting  for $G$
if $\lambda (G, \vec{y})=\lambda(G)$.
The following fact is easily implied by the definition of the Lagrangian.

\begin{fact}\label{mono}
Let $G_1$, $G_2$ be $r$-uniform graphs and $G_1\subseteq G_2$. Then $\lambda (G_1) \le \lambda (G_2).$
\end{fact}
In \cite{MS}, Motzkin and Straus established a remarkable connection between the clique number and the Lagrangian of a graph.

\begin{theo} \cite{MS}  \label{MStheo}
If $G$ is a 2-graph in which a maximum clique has order $t$ then
$\lambda(G)=\lambda(K^{(2)}_t)={1 \over 2}(1 - {1 \over t})$.
\end{theo}

%This new proof aroused interests in the study of Lagrangians of hypergraphs.
The Motzkin-Straus result provides solutions to the optimization problem of a class of homogeneous  multilinear functions  over the standard simplex of the  Euclidean  space. The Motzkin-Straus result
and its extension were also successfully employed in optimization to provide heuristics for the
maximum clique problem (see \cite{B1,B2,B3,G9,PP15}). It is interesting to explore whether  similar results holds for hypergraphs. The obvious generalization of Motzkin and Straus' result to hypergraphs is false because there are many examples of hypergraphs that do not achieve their
Lagrangian on any proper subhypergraph.

Lagrangians of hypergraphs has been proved to be a useful tool in hypergraph extremal problems. Applications of Lagrangian method can be found in \cite{FF,FR84,keevash,mubayi06,sidorenko89}. In most applications, an upper bound is needed. Frankl and F\"uredi \cite{FF} asked the following question. Given $r \ge 3$ and $m \in {\mathbb N}$ how large can the Lagrangian of an $r$-graph with $m$ edges be?
For distinct $A, B \in {\mathbb N}^{(r)}$ we say that $A$ is less than $B$ in the {\em colex ordering} if $max(A \triangle B) \in B$, where $A \triangle B=(A \setminus B)\cup (B \setminus A)$. For example, the first $t \choose r$ $r$-tuples in the colex ordering of ${\mathbb N}^{(r)}$ are the edges of $[t]^{(r)}$. The following conjecture of Frankl and F\"uredi (if it is true) proposes a  solution to the question mentioned above.

\begin{con}  \cite{FF} \label{conjecture} The $r$-graph with $m$ edges formed by taking the first $m$ sets in the colex ordering of ${\mathbb N}^{(r)}$ has the largest Lagrangian of all $r$-graphs with  $m$ edges. In particular, the $r$-graph with $t \choose r$ edges and the largest Lagrangian is $[t]^{(r)}$.
\end{con}

This conjecture is true when $r=2$ by Theorem \ref{MStheo}. For the case $r=3$, Talbot in \cite{T} proved the following.

\begin{theo}  \cite{T}  \label{Tal} Let $m$ and $t$ be integers satisfying
${t-1 \choose 3} \le m \le {t-1 \choose 3} + {t-2 \choose 2} - (t-1).$
Then Conjecture \ref{conjecture} is true for $r=3$ and this value of $m$.
%Conjecture \ref{conjecture} is also true for $r=3$ and $m= {t \choose 3}-1$ or $m={t \choose 3} -2$.
\end{theo}

Recently, in \cite{TPZZ}, using some different approaches,  Conjecture \ref{conjecture} is confirmed for $r=3$ when the value of $m$ satisfying
${t-1 \choose 3} \le m \le {t-1 \choose 3} + {t-2 \choose 2} - \frac{1}{2}(t-1).$

%In \cite{T}, the following result is also proved, which is the evidence for Conjecture \ref{conjecture} for
%$r$-graphs $G$ on exactly $t$ vertices.

%\begin{theo} (Talbot \cite{T}) \label{Talr} For any $r \ge 4$ there exists constants $\gamma_r$ and $\kappa_0(r)$ such that if $m$ satisfies
%$${t-1 \choose r} \le m \le {t-1 \choose r} + {t-2 \choose r-1} - \gamma_r (t-1)^{r-2},$$
%with $t \ge \kappa_0(r)$, let $G$ be an $r$-graph on $t$ vertices with $m$ edges, then $\lambda(G)\leq\lambda([t-1]^{(r)})$.
%\end{theo}

%The above result  improves Theorem \ref{Tal}.

Although the obvious generalization of Motzkin and Straus' result to hypergraphs is false as mentioned earlier,  we attempt to explore the relationship between the Lagrangian of a hypergraph and the  size of its maximum cliques for hypergraphs when the number of edges is in certain ranges. In \cite{PZ}, it is conjectured that the following Motzkin and Straus type results are true for hypergraphs.

\begin{con} \label{conjecture1} \cite{PZ}
Let $t$, $m$, and $r\ge 3$ be positive integers satisfying ${t-1 \choose r} \le m \le {t-1 \choose r} + {t-2 \choose r-1}$.
Let $G$ be an $r$-graph with $m$ edges and $G$ contain a clique of order  $t-1$. Then $\lambda(G)=\lambda([t-1]^{(r)})$.
\end{con}

\begin{con} \label{conjecture2} \cite{PZ} Let $t$, $m$, and $r\ge 3$ be positive integers satisfying ${t-1 \choose r} \le m \le {t-1 \choose r} + {t-2 \choose r-1}$.
Let $G$ be an $r$-graph with $m$ edges without containing a clique of order  $t-1$. Then $\lambda(G) < \lambda([t-1]^{(r)})$.
\end{con}

Note that the upper bound ${t-1 \choose r} + {t-2 \choose r-1}$ in Conjecture \ref{conjecture1} is the best possible (see \cite{PZ}).  Conjecture \ref{conjecture1} is confirmed when $r=3$ in \cite{PZ}. Let $C_{r,m}$ denote the $r$-graph with $m$ edges formed by taking the first $m$ sets in the colex ordering of ${\mathbb N}^{(r)}$.
The following result was given in  \cite{T}.

\begin{lemma} \cite{T} \label{LemmaTal7}
For any integers $m, t,$ and $r$ satisfying ${t-1 \choose r} \le m \le {t-1 \choose r} + {t-2 \choose r-1}$,
we have $\lambda(C_{r,m}) = \lambda([t-1]^{(r)})$.
\end{lemma}

In \cite{PTZ}, the following result is obtained for $r$-graphs.

\begin{theo} \cite{PTZ}\label{TheoremPTZ}
 Let $t$,$m$ and $r$ be positive integers satisfying ${t-1 \choose r} \le m \le {t-1 \choose r} + {t-2 \choose r-1}-(2^{r-3}-1)({t-2 \choose r-2}-1)$.  Let $G$ be an $r$-graph with $t$ vertices and $m$ edges and contain a clique of order $t-1$. Then $\lambda(G) = \lambda([t-1]^{(r)})$.
\end{theo}

In \cite{T}, the following result is also proved, which is the evidence for Conjecture \ref{conjecture} for
$r$-graphs $G$ on exactly $t$ vertices.

\begin{theo}  \cite{T}  \label{Talr} For any $r \ge 4$ there exists constants $\gamma_r$ and $\kappa_0(r)$ such that if $m$ satisfies
$${t-1 \choose r} \le m \le {t-1 \choose r} + {t-2 \choose r-1} - \gamma_r (t-1)^{r-2},$$
with $t \ge \kappa_0(r)$, let $G$ be an $r$-graph on $t$ vertices with $m$ edges, then $\lambda(G)\leq\lambda([t-1]^{(r)})$.
\end{theo}
The main result in this paper is Theorem \ref{mainresult} which is a accompany result of Theorem \ref{TheoremPTZ}.

\begin{theo} \label{mainresult} Let $m$, $t$, and $r\geq 4$ be integers satisfying
${t-1\choose r}\leq m \leq {t-1\choose r}+ {t-2\choose r-1}-[(2r-6)\times2^{r-1}+2^{r-3}+(r-4)(2r-7)-1]({t-2\choose r-2}-1).$ Let $G$ be an $r$-graph with $t$ vertices and $m$ edges  and without containing a clique of order $t-1$. Then $\lambda(G) < \lambda([t-1]^{(r)})$.
\end{theo}
Theorem \ref{mainresult} and Theorem \ref{TheoremPTZ}  give a Motzkin-Straus result for some $r$-graph. Combing Theorems \ref{TheoremPTZ}  and \ref{mainresult}, we have the following result immediately.
\begin{coro} \label{improvement} Let $m$, $t$, and $r\geq 4$  be integers satisfying
${t-1\choose r}\leq m \leq {t-1\choose r}+ {t-2\choose r-1}-[(2r-6)\times2^{r-1}+2^{r-3}+(r-4)(2r-7)-1]({t-2\choose r-2}-1).$ Let $G$ be an $r$-graph with $t$ vertices and  $m$ edges. Then $\lambda(G) \leq \lambda([t-1]^{(r)})$.
\end{coro}
Note that ${t-1\choose r}\leq m \leq {t-1\choose r}+ {t-2\choose r-1}-[(2r-6)\times2^{r-1}+2^{r-3}+(r-4)(2r-7)-1]({t-2\choose r-2}-1)$ implies   the number of vertices $t$ should be sufficiently large such that  ${t-2\choose r-1}\geq[(2r-6)\times2^{r-1}+2^{r-3}+(r-4)(2r-7)-1]({t-2\choose r-2}-1)$ in Theorem \ref{mainresult} and Corollary \ref{improvement}.
%Combing Theorems \ref{TheoremPTZ} and \ref{mainresult}, we have the following result immediately.
%For Conjecture \ref{conjecture1}, we have

Theorem \ref{mainresult} and  Corollary \ref{improvement} provide evidence for both Conjecture \ref{conjecture2} and Conjecture \ref{conjecture} respectively. The contribution of Corollary \ref{improvement} is that the method developed in the proof of  Theorem \ref{mainresult} is simpler and different from that in Theorem \ref{Talr} in some ways. The upper bound in Corollary \ref{improvement} for the number of edges $m$ is more explicit and an improvement comparing to the bound in Theorem \ref{Talr}.  The proof of  Theorem \ref{mainresult}  will be given in Section \ref{section1}. Further remarks and conclusions are in Section \ref{section3}.
%Next we give some useful results.

%\section{Useful Results}

%In other words, for any $i<j$, if $k_1 k_2\ldots k_{r-1}  \in E_j$, where $k_1, k_2, \ldots, k_{r-1} \neq i$, then  $k_1k_2 \ldots k_{r-1}  \in E_i$.

\section{Proof of Theorem \ref{mainresult}}\label{section1}

We will impose one additional condition on any optimal weighting ${\vec x}=(x_1, x_2, \ldots, x_n)$ for an $r$-graph $G$:
\begin{eqnarray}
 &&|\{i : x_i > 0 \}|{\rm \ is \ minimal, i.e. \ if}  \ \vec y {\rm \ is \ a \ feasible \ weighting \ for \ } G  {\rm \ satisfying }\nonumber \\
 &&|\{i : y_i > 0 \}| < |\{i : x_i > 0 \}|,  {\rm \  then \ } \lambda (G, {\vec y}) < \lambda(G) \label{conditionb}.
\end{eqnarray}

When the theory of Lagrange multipliers is applied to find the optimum of $\lambda(G, {\vec x})$, subject to $\sum_{i=1}^n x_i =1$, notice that $\lambda (E_i, {\vec x})$ corresponds to the partial derivative of  $\lambda(G, \vec x)$ with respect to $x_i$. The following lemma gives some necessary conditions of an optimal weighting for $G$.

\begin{lemma}   \cite{FR84}   \label{LemmaTal5} Let $G=(V,E)$ be an $r$-graph on the vertex set $[n]$ and ${\vec x}=(x_1, x_2, \ldots, x_n)$ be an optimal  weighting for $G$ with $k$  ($\le n$) non-zero weights $x_1$, $x_2$, $\cdots$, $x_k$  satisfying condition (\ref{conditionb}). Then for every $\{i, j\} \in [k]^{(2)}$, (a) $\lambda (E_i, {\vec x})=\lambda (E_j, \vec{x})=r\lambda(G)$, (b) there is an edge in $E$ containing both $i$ and $j$.
\end{lemma}

\begin{defi}
An $r$-graph $G=(V,E)$ on the vertex set $[n]$ is left-compressed if $j_1j_2\ldots j_r\in E$ implies $i_1i_2\ldots i_r\in E$ whenever $i_k \leq j_k, 1\leq k \leq r$. Equivalently, an $r$-graph $G=(V,E)$  on the vertex set $[n]$ is  left-compressed if $E_{j\setminus i}=\emptyset$ for any $1\le i<j\le n$.
\end{defi}

\begin{remark}\label{r1} (a) In Lemma \ref{LemmaTal5}, part (a) implies that
$x_j\lambda(E_{ij}, {\vec x})+\lambda (E_{i\setminus j}, {\vec x})=x_i\lambda(E_{ij}, {\vec x})+\lambda (E_{j\setminus i}, {\vec x}).$
In particular, if $G$ is left-compressed, then
$(x_i-x_j)\lambda(E_{ij}, {\vec x})=\lambda (E_{i\setminus j}, {\vec x})$
for any $i, j$ satisfying $1\le i<j\le k$ since $E_{j\setminus i}=\emptyset$.

(b) If  $G$ is left-compressed, then for any $i, j$ satisfying $1\le i<j\le k$,
\begin{equation}\label{enbhd}
x_i-x_j={\lambda (E_{i\setminus j}, {\vec x}) \over \lambda(E_{ij}, {\vec x})}
\end{equation}
holds.  If  $G$ is left-compressed and  $E_{i\setminus j}=\emptyset$ for $i, j$ satisfying $1\le i<j\le k$, then $x_i=x_j$.

(c) By (\ref{enbhd}), if  $G$ is left-compressed, then an optimal  weighting  ${\vec x}=(x_1, x_2, \ldots, x_n)$ for $G$  must satisfy
%\begin{equation}\label{conditiona}
$x_1 \ge x_2 \ge \ldots \ge x_n \ge 0$.
%\end{equation}
\end{remark}
%Denote $$\lambda_{m}^{r}=\max\{\lambda(G): G {\rm \ is \ an \ } r-{\rm graph\ with \ } m {\rm \ edges }\}.$$
 %We will use Lemmas  \ref{LemmaTal8} and \ref{TheoPTZ} in proving our main results.
%\begin{lemma} (Peng, Tang, Zhao \cite{PTZ}) \label{TheoPTZ} Let $m$ and $t$ be positive integers satisfying ${t-1 \choose 3} \le m \le {t-1 \choose 3} + {t-2 \choose 2}$. Let $G=(V,E)$ be a left-compressed $3$-graph with $m$ edges and $t$ vertices satisfying
%$|[t-2r+6]^{(2)}\setminus E_{t}| \ge |E_{(t-1)t}|$.
%Then $\lambda(G) \leq \lambda([t-1]^{(r-1)})$.
%\end{lemma}
%The following lemma implies that we only need to consider left-compressed $r$-graphs when Conjecture \ref{conjecture} is explored.
%\begin{lemma} (Talbot \cite{T}) \label{LemmaTal8}
%Let $m,t$ be positive integers satisfying  $m\leq{t \choose r}-1$, then there exists a left-compressed $r$-graph $G$ with $m$ edges such that $\lambda(G)=\lambda_{m}^{r}.$
%\end{lemma}
%Denote $\lambda_{m}^{r}=\max \{ \lambda(G):$
%$ G$ is an  $r$-graph    with $m$   edges$\}$.
%\begin{lemma}\cite{T}\label{leftcom0}
%There exists a left-compressed $r$-graph $G$  with $m$ edges   such that
%$\lambda(G)=\lambda_{m}^{r}$.
%\end{lemma}

%\section{Proof of Theorem \ref{mainresult}}\label{section1}

%In this section, we prove that Conjecture \ref{conjecture} and Conjecture \ref{conjecture2} hold when  $m$ lies in a certain range.ma
Denote $\lambda_{(m,  t)}^{r}=\max \{ \lambda(G):$
$ G$ is an  $r$-graph with $t$ vertices  and $m$   edges $  \}$. The following Lemma is proved in \cite{T}.

\begin{lemma}\cite{T}\label{leftcom0}
There exists a left-compressed $r$-graph $G$  with $t$ vertices and $m$ edges   such that
$\lambda(G)=\lambda_{(m,t)}^{r}$.
\end{lemma}
\begin{remark}\label{r2}
Since the only left-compressed $r$-graph with $t$ vertices and $m={t\choose r}$ edges is $[t]^{(r)}$. Hence by Lemma \ref{leftcom0} and Fact \ref{mono}, we have $\lambda_{(m,t)}^{r}\leq\lambda([t]^{(r)})$.
\end{remark}
Denote $\lambda_{(m, t-1, t)}^{r-}=\max \{ \lambda(G):$
$ G$ is an  $r$-graph with $t$ vertices  and $m$   edges    not containing a clique of order $t-1  \}$.  The following lemma  implies that we only need to consider left-compressed  $r$-graphs $G$ when we prove Theorem \ref{mainresult}.

\begin{lemma}\label{leftcom} Let $m$ and $t$ be integers satisfying
${t-1\choose r}\leq m \leq {t-1\choose r}+ {t-2\choose r-1}-[(2r-6)\times2^{r-1}+2^{r-3}+(r-4)(2r-7)-1]({t-2\choose r-2}-1).$
There exists a left-compressed $r$-graph $G$ on vertex set $[t]$ with $m$ edges   without containing  $[t-1]^{(r)}$ such that
$\lambda(G)=\lambda_{(m, t-1, t)}^{r-}$.
\end{lemma}
In the proof of Lemma \ref{leftcom}, we need to define some partial order relation.
An $r$-tuple  $i_1 i_2\cdots i_r$ is called a   descendant   of an $r$-tuple  $j_1j_2\cdots j_r$ if $i_s\le j_s$ for each $1\le s\le r$, and $i_1+i_2+\cdots +i_r < j_1+j_2+\cdots +j_r$. In this case, the $r$-tuple $j_1j_2\cdots j_r$   is called an   ancestor  of $i_1 i_2\cdots i_r$.  The $r$-tuple $i_1i_2\cdots i_r$   is called a  direct descendant  of $j_1 j_2\cdots j_r$ if   $i_1i_2\cdots i_r$   is a  descendant of $j_1j_2\cdots j_r$ and $j_1+j_2+\cdots +j_r=i_1+i_2+\cdots +i_r +1$.  We say that $i_1 i_2\cdots i_r$ has lower hierarchy than  $j_1j_2\cdots j_r$ if $i_1 i_2\cdots i_r$  is  a descendant of $j_1j_2\cdots j_r$. This is a partial order on the set of all $r$-tuples.

\noindent {\em Proof of Lemma \ref{leftcom}.} Let $G$ be an $r$-graph with $t$ vertices  and $m$  edges  without containing a clique of order $t-1$  such that $\lambda(G)=\lambda_{(m, t-1, t)}^{r-}$. We call  $G$ an extremal $r$-graph for $m$, $t-1$ and $t$. Let ${\vec x}=(x_1, x_2, \ldots, x_t)$ be an optimal weighting of $G$.    We can assume that $x_i\ge x_j$ when $i<j$ since otherwise we can just relabel the vertices of $G$ and obtain another  extremal $r$-graph for $m$,  $t-1$ and $t$ with an optimal weighting ${\vec x}=(x_1, x_2, \ldots, x_t)$ satisfying $x_i\ge x_j$ when $i<j$. Next we obtain a new $r$-graph $H$ from $G$ by performing the following:

\begin{enumerate}

\item If $(t-r)\ldots(t-1) \in E(G)$, then there is at least one $r$-tuple  in $[t-1]^{(r)}\setminus E(G)$,  we replacing $(t-r)\ldots(t-1)$ by this $r$-tuple;

\item If  an  edge in $G$ has  a  descendant other than $(t-r)\ldots(t-1)$  that is not  in $E(G)$, then replace this edge by a descendant other than $(t-r)\ldots(t-1)$ with the lowest hierarchy. Repeat this until there is no such an edge.

\end{enumerate}

  Then $H$ satisfies the following properties:

\begin{enumerate}

\item  The number of edges in $H$ is the same as the number of edges in $G$.

\item $\lambda(G)=\lambda(G, {\vec x})\le \lambda(H, {\vec x})\le \lambda(H).$

\item $(t-r)\ldots(t-1) \notin E(H)$.

\item For any edge in $E(H)$, all its  descendants  other than $(t-r)\ldots(t-1)$ will be in $E(H)$.

\end{enumerate}

If $H$ is not left-compressed, then there is an ancestor  of $(t-r)\ldots(t-1)$, says $e$, such that $e\in E(H)$. Hence $(t-r)\ldots(t-2)t$ and all the descendants of $(t-r)\ldots(t-2)t$ other than $(t-r)\ldots(t-1)$ will be in $E(H)$. Then
$$m\ge {t-1 \choose r}-1 + {t-2 \choose r-1}>{t-1\choose r}+ {t-2\choose r-1}-[(2r-6)\times2^{r-1}+2^{r-3}+(r-4)(2r-7)-1]({t-2\choose r-2}-1)$$ which is a contradiction.   $H$ does not contain $[t-1]^{(r)}$ since $H$ does not contain $(t-r)\ldots(t-1)$. Clearly  $H$ is on vertex set $[t]$. So we complete the proof of Lemma \ref{leftcom}\qed

In the rest of the paper we assume that $r\geq 4$ be an integer. In the following three lemmas, Lemma \ref{Lemmaeq} implies the  maximum weight of $G$ should distribute  'uniform' on the $t$ vertices if $\lambda(G)\geq \lambda([t-1]^{(r)})$, and  Lemma \ref{lemmac} implies $G$ contains most of the first ${t-2r+6 \choose r}$ edges in colex ordering of $N^{(r)}$ if $\lambda(G) \geq \lambda([t-1]^{(r)}),$  while Lemma \ref{lemmab} implies $G$ also contains most of  the first ${t-2r+6 \choose r-1}$ edges containing $t-1$. Since $G$ is left-compressed, $G$ also contains most of the the first ${t-2r+6 \choose r-1}$ edges containing vertex $i$, where $t-2r+7\leq i\leq t-1.$ So $G$ contains most edges of $[t-1]^{(r)}$. Note that,  in the proof of Lemma 2.6, whenever the lower bound of a product is greater than the upper bound, we take this to be the empty product.

\begin{lemma}\label{Lemmaeq}  (a) Let $G$ be an $r$-graph on vertex set $[t]$. Let  $\vec{x}=(x_{1},x_{2},\ldots ,x_{t})$ be an optimal weighting for $G$ satisfying $x_1 \ge x_2 \ge \ldots \ge x_t \ge 0$. Then
$x_1< x_{t-2r+3}+x_{t-2r+4}$  or  $$\lambda(G)\le \frac{1}{r!}\frac{(t-r)^{r-1} \prod\limits_{i=t-r+2}^{t-2}i}{(t-r+1)^{r-2}(t-1)^{r-2}}<\frac{1}{r!}\frac{ \prod\limits_{i=t-r}^{t-1}i}{(t-1)^r}=\lambda([t-1]^{(r)}).$$

(b)  Let $G$ be an $r$-graph on vertex set $[t]$. Let  $\vec{x}=(x_{1},x_{2},\ldots ,x_{t})$ be an optimal weighting for $G$ satisfying $x_1 \ge x_2 \ge \ldots \ge x_t \ge 0$. Then
$x_1<2( x_{t-2r+4}+x_{t-2r+5})$  or
$$\lambda(G)\le \frac{1}{r!}\frac{(t-r)^{r-1} \prod\limits_{i=t-r+2}^{t-2}i}{(t-r+1)^{r-2}(t-1)^{r-2}}<\frac{1}{r!}\frac{ \prod\limits_{i=t-r}^{t-1}i}{(t-1)^r}=\lambda([t-1]^{(r)}).$$

\end{lemma}
{\em Proof.} (a) If $x_1\geq x_{t-2r+3}+x_{t-2r+4}$,  then $rx_1+x_2+\cdots+x_{t-2r+2}\geq x_1+x_2+\cdots +x_{t-2r+4}+x_{t-3}+x_{t-2r+6}+x_{t-1}+x_{t}=1$. Recalling that $x_1 \ge x_2 \ge \ldots \ge x_{t-2r+2}$, we have $x_1\geq \frac{1}{t-r+1}$.
Using Lemma \ref{LemmaTal5}, we have $\lambda(G)=\frac{1}{r}\lambda(E_1,x).$ Note that $E_1$ is an  $(r-1)$-graph with $t-1$ vertices and total weights at most $1-\frac{1}{t-r+1}$. Hence by Remark \ref{r2}(change the total weights 1 to $1-\frac{1}{t-r+1}$).
\begin{eqnarray}\label{eq234}
\lambda(G)&=&\frac{1}{r}\lambda(E_1,x)\leq\frac{1}{r}{t-1 \choose r-1}(\frac{1-\frac{1}{t-r+1}}{t-1})^{r-1}\nonumber\\
&=&\frac{1}{r!}\frac{(t-r)^{r-1} \prod\limits_{i=t-r+2}^{t-2}i}{(t-r+1)^{r-2}(t-1)^{r-2}}.
\end{eqnarray}
Next we prove
\begin{eqnarray}\label{eqa}
\frac{1}{r!}\frac{(t-r)^{r-1} \prod\limits_{i=t-r+2}^{t-2}i}{(t-r+1)^{r-2}(t-1)^{r-2}}
< \frac{1}{r!}\frac{ \prod\limits_{i=t-r}^{t-1}i}{(t-1)^r}
=\lambda([t-1]^{(r)}).
\end{eqnarray}
To show this, we only need to prove
\begin{eqnarray}\label{eqb}
(t-r)^{r-2}(t-1)<(t-r+1)^{r-1}.
\end{eqnarray}
If $t=r,r+1$, (\ref{eqb}) clearly holds. Assuming $t\geq r+2$,  we prove this inequality by induction.
Now we suppose that (\ref{eqb}) holds for some $r\geq 4$, we will show it also holds for $r+1$.
Replacing $t$ by $t-1$ in (\ref{eqb}). We have
$$[t-(r+1)]^{r-2}(t-2)<(t-r)^{r-1}.$$
Multiplying $t-(r+1)$ to the above inequality, we have
$$[t-(r+1)]^{r-1}(t-2)<(t-r)^{r-1}[t-(r+1)].$$
Adding $[t-(r+1)]^{r-1}$ to the above inequality, we obtain
\begin{eqnarray}
[t-(r+1)]^{r-1}(t-1)&<&(t-r)^{r-1}[t-(r+1)]+[t-(r+1)]^{r-1}\nonumber\\
&=&(t-r)^{r}-(t-r)^{r-1}+[t-(r+1)]^{r-1}<(t-r)^{r}.
\end{eqnarray}
Hence (\ref{eqb}) also holds for $r+1$ and the induction is complete.

(b) If $x_1\geq 2( x_{t-2r+5}+x_{t-2r+6})$,  then $x_1+x_2+\cdots+x_{t-2r+4}+(r-2)\frac{x_1}{2}\geq x_1+x_2+\cdots +x_{t-2r+4}+x_{t-3}+x_{t-2r+6}+x_{t-1}+x_{t}=1$. Recalling that $x_1 \ge x_2 \ge \ldots \ge x_{t-2r+4}$ and $r\geq 4$, we have $x_1\geq \frac{1}{t-2r+4+\frac{r-2}{2}}\geq \frac{1}{t-r+1}$. The rest of the proof is  identical to that in part (a), we omit the computation details here.
\qed

\begin{lemma} \label{lemmab}
% satisfying ${t-1 \choose 3} \le m \le {t-1 \choose 3} + {t-2 \choose 2}$.
Let $G$ be a left-compressed $r$-graph  on the vertex set $[t]$ without containing $[t-1]^{(r)}$, then $|[t-2r+6]^{(r-1)} \backslash E_{t-1}| \leq  2^{r-1}|E_{(t-1)t}|$ or $\lambda(G) < \lambda([t-1]^{(r)}).$
 \end{lemma}
\noindent{\em Proof.}
Let  ${\vec x}=(x_1, x_2, \ldots, x_t)$ be an optimal weighting for $G$. Since $G$ is left-compressed, by Remark \ref{r1}(a), $x_1\ge x_2 \ge \cdots \ge x_t \ge 0$. If $x_{t}=0$, then $\lambda(G)=\lambda(G, \vec{x}) < \lambda([t-1]^{(r)})$ since $G$ does not contain $[t-1]^{(r)}.$ So we assume that $x_{t}>0$.

%Now we proceed to show that $\lambda(G) \le \lambda([t-1]^{(r)})$.
Consider a new weighting for $G$, ${\vec y}=(y_1, y_2, \ldots, y_t)$ given by $y_i=x_i$ for $i\neq t-1, t$, $y_{t-1}=x_{t-1}+x_t$ and $y_t=0$. By Lemma \ref{LemmaTal5}(a), $\lambda(E_{t-1}, \vec{x})=\lambda(E_{t}, \vec{x})$, so
\begin{eqnarray}\label{eq10b}
\lambda(G,\vec {y})- \lambda(G,\vec {x})&=&x_{t}(\lambda(E_{t-1}, \vec{x})-x_{t}\lambda(E_{(t-1)t}, \vec{x})) \nonumber \\
&&-x_{t}(\lambda(E_{t}, \vec{x})-x_{t}\lambda(E_{(t-1)t}, \vec{x}))-x_{t-1}x_t\lambda(E_{(t-1)t}, \vec{x}))\nonumber \\
&=&x_{t}(\lambda(E_{t-1}, \vec{x})-\lambda(E_{t}, \vec{x}))-x_{t}^2\lambda(E_{(t-1)t}, \vec{x})  \nonumber\\
&=& -x_{t}^2\lambda(E_{(t-1)t}, \vec{x}).
\end{eqnarray}

%Since $y_{t-1}=0$ we may remove all edges containing $t-1$ from $E$ to form a new $r$-graph $G^*=([t], E^*)$ with $\lambda(G^*,\vec {y})=
%\lambda(G,\vec {y})$ and $\vert E^*\vert=\vert E\vert-\vert E_{t-1}\vert$.

Assume that $|[t-2r+6]^{(r-1)} \backslash E_{t-1}| >  2^{r-1}|E_{(t-1)t}|$. If $\lambda(G) < \lambda([t-1]^{(r)})$ we are done. Otherwise if $\lambda(G) \geq \lambda([t-1]^{(r)})$ we will show that there exists a set of edges $F\subset [t-1]^{(r)}\setminus E$ satisfying
\begin{equation}\label{eq11b}
\lambda(F,\vec {y})> x_{t}^2\lambda(E_{(t-1)t}, \vec{x}).
\end{equation}

Then using (\ref{eq10b}) and (\ref{eq11b}), the $r$-graph $G^{*}=([t], E^{*})$, where $E^{*}=E\cup F$, satisfies $\lambda(G^{*}, \vec {y})=\lambda(G, \vec {y})+\lambda(F, \vec {y})> \lambda(G, \vec{x})=\lambda(G)$. Since $\vec {y}$ has only $t-1$ positive weights, then $\lambda(G^{*}, \vec {y})\le \lambda([t-1]^{(r)})$, and consequently,
$\lambda(G)<\lambda([t-1]^{(r)}).$ This is a contradiction.

We now construct the set of edges $F$. Let $C=[t-2r+6]^{(r-1)} \setminus E_{t-1}$. Then by the assumption,
$\vert C\vert > 2^{r-1}|E_{(t-1)t}|$ and
$\lambda(C, \vec{x})\ge 2^{r-1}|E_{(t-1)t}|x_{t-3r+8}\ldots x_{t-2r+6}.$

Let $F$ consist of those edges in $[t-1]^{(r)}\setminus E$ containing the vertex $t-1$.  Since $\lambda(G) \geq \lambda([t-1]^{(r)})$ then $x_{t-2r+3}> \frac{x_1}{2}$ by Lemma \ref{Lemmaeq}(a) and $x_{t-2r+4}\geq x_{t-2r+5} >\frac{x_1}{4}$ by Lemma \ref{Lemmaeq}(b). Hence
\begin{eqnarray}
\lambda(F,\vec {y}) &=&(x_{t-1}+x_{t})\lambda(C, \vec{x}) > 2x_{t}\cdot 2^{r-1}|E_{(t-1)t}|x_{t-3r+8}\ldots x_{t-2r+6} \nonumber \\
&\ge&  x_{t}^2 |E_{(t-1)t}|(x_1)^{2} \ge x_{t}^2 \sum_{i_1\ldots i_{r-2}\in E_{(t-1)t}} x_{i_1}\ldots x_{i_2} = x_{t}^2 \lambda(E_{(t-1)t}, \vec{x}).
\end{eqnarray}
Hence $F$ satisfies (\ref{eq11b}). This proves Lemma \ref{lemmab}.\qed

\begin{lemma} \label{lemmac}
% satisfying ${t-1 \choose 3} \le m \le {t-1 \choose 3} + {t-2 \choose 2}$.
Let $G$ be a left-compressed $r$-graph  on the vertex set $[t]$ without containing $[t-1]^{(r)}$, then $|[t-2r+6]^{(r)} \backslash E| \leq  2^{r-1}|E_{(t-1)t}|$ or $\lambda(G) < \lambda([t-1]^{(r)}).$
 \end{lemma}
\noindent{\em Proof.}
%The proof of this lemma is similar to the proof of Lemma \ref{lemmab}.
Let  ${\vec x}=(x_1, x_2, \ldots, x_t)$ be an optimal weighting for $G$. Since $G$ is left-compressed, by Remark \ref{r1}(a), $x_1\ge x_2 \ge \cdots \ge x_t \ge 0$. If $x_{t}=0$, then $\lambda(G) < \lambda([t-1]^{(r)})$ since $G$ does not contain  $[t-1]^{(r)}.$ So we assume that $x_{t}>0$.

%Now we proceed to show that $\lambda(G) \le \lambda([t-1]^{(r)})$.
Consider a new weighting for $G$, ${\vec y}=(y_1, y_2, \ldots, y_t)$ given by $y_i=x_i$ for $i\neq t-1, t$, $y_{t-1}=x_{t-1}+x_t$ and $y_t=0$. By Lemma \ref{LemmaTal5}(a), $\lambda(E_{t-1}, \vec{x})=\lambda(E_{t}, \vec{x})$, similar to (4), we have
\begin{eqnarray}\label{eq10c}
\lambda(G,\vec {y})- \lambda(G,\vec {x})=-x_{t}^2\lambda(E_{(t-1)t}, \vec{x}).
\end{eqnarray}

%Since $y_{t-1}=0$ we may remove all edges containing $t-1$ from $E$ to form a new $r$-graph $G^*=([t], E^*)$ with $\lambda(G^*,\vec {y})=
%\lambda(G,\vec {y})$ and $\vert E^*\vert=\vert E\vert-\vert E_{t-1}\vert$.

Assume that $|[t-2r+6]^{(r)} \backslash E| >  2^{r-1}|E_{(t-1)t}|$. If $\lambda(G) < \lambda([t-1]^{(r)})$ we are done. Otherwise if $\lambda(G) \geq \lambda([t-1]^{(r)})$ we will show that there exists a set of edges $F\subset [t-2r+6]^{(4)}\setminus E$ satisfying
\begin{equation}\label{eq11c}
\lambda(F,\vec {y})> x_{t}^2\lambda(E_{(t-1)t}, \vec{x}).
\end{equation}

Then using (\ref{eq10c}) and (\ref{eq11c}), the $r$-graph $G^{*}=([t], E^{*})$, where $E^{*}=E\cup F$, satisfies $\lambda(G^{*}, \vec {y})=\lambda(G, \vec {y})+\lambda(F, \vec {y})> \lambda(G, \vec{x})=\lambda(G)$. Since $\vec {y}$ has only $t-1$ positive weights, then $\lambda(G^{*}, \vec {y})\le \lambda([t-1]^{(r)})$, and consequently,
$\lambda(G)<\lambda([t-1]^{(r)}).$ This is a contradiction.

We now construct the set of edges $F$. Let $C=[t-2r+6]^{(r)} \setminus E$. Then by the assumption,
$\vert C\vert > 2^{r-1}|E_{(t-1)t}|$ and
$\lambda(C, \vec{x})\ge 2^{r-1}|E_{(t-1)t}|x_{t-3r+7}\ldots x_{t-2r+6}.$

Let $F=C$.  Since $\lambda(G) \geq \lambda([t-1]^{(r)})$ then $x_{t-2r+3}\geq \frac{x_1}{2}$ by Lemma \ref{Lemmaeq}(a) and $x_{t-2r+4}\geq x_{t-2r+5}>\frac{x_1}{4}$ by Lemma \ref{Lemmaeq}(b). Hence
\begin{eqnarray}
\lambda(F,\vec {y}) &=&\lambda(C, \vec{x}) > 2^{r-1}|E_{(t-1)t}|x_{t-3r+7}\ldots x_{t-2r+6} \ge  x_{t}^2 |E_{(t-1)t}|(x_1)^{2}   \nonumber \\
&\ge& x_{t}^2 \sum_{i_1\ldots i_{r-2}\in E_{(t-1)t}} x_{i_1}\ldots x_{i_{r-2}} = x_{t}^2 \lambda(E_{(t-1)t}, \vec{x}).
\end{eqnarray}
Hence $F$ satisfies (\ref{eq11c}). This proves Lemma \ref{lemmac}.\qed

Now we are ready to  prove Theorem \ref{mainresult}.

\noindent{\em Proof of Theorem \ref{mainresult}.}
Let $m$ and $t$ be integers satisfying
${t-1\choose r}\leq m \leq {t-1\choose r}+ {t-2\choose r-1}-[(2r-6)\times2^{r-1}+2^{r-3}+(r-4)(2r-7)-1]({t-2\choose r-2}-1).$ Let $G$ be an $r$-graph  with $t$ vertices and $m$ edges without containing a clique of order $t-1$ such that $\lambda(G)=\lambda_{(m, t-1, t)}^{r-}$. Then by Lemma \ref{leftcom}, we can assume that $G$ is left-compressed and does not contain $[t-1]^{(r)}$. Let  ${\vec x}=(x_1, x_2, \ldots, x_t)$ be an optimal weighting for $G$. Since $G$ is left-compressed, by Remark \ref{r1}(a), $x_1\ge x_2 \ge \cdots \ge x_t \ge 0$. If $x_{t}=0$, then $\lambda(G) < \lambda([t-1]^{(r)})$ since $G$ does not contain $[t-1]^{(r)}$. So we assume that $x_{t}>0$.

If $\lambda(G) < \lambda([t-1]^{(r)})$ we are done. Otherwise $|[t-2r+6]^{(r-1)} \backslash E_{t-1}| \leq  2^{r-1}|E_{(t-1)t}|$ by Lemma \ref{lemmab}. Recalling that $G$ is left-compressed, we have $|[t-2r+6]^{(r-1)} \backslash E_{i}| \leq  2^{r-1}|E_{(t-1)t}|$ for $t-2r+7\leq i\leq t-1.$ We also have $|[t-2r+6]^{(4)} \backslash E| \leq 2^{r-1}|E_{(t-1)t}|$ by Lemma \ref{lemmac}. Note that  $|E_{(t-1)t}|\leq{t-2\choose r-2}-1$,  then
\begin{eqnarray}
|[t-1]^{(r)}\bigcap E|&\geq &|[t-2r+6]^{(r)} \bigcap E|+\sum\limits_{i=t-2r+7}^{t-1}|[t-2r+6]^{(r-1)} \bigcap E_{i}|\nonumber\\
&\geq & {t-2r+6 \choose r}-2^{r-1}|E_{(t-1)t}|+(2r-7)({t-2r+6 \choose r-1}-(2r-7)\times 2^{r-1}|E_{(t-1)t}|)\nonumber\\
&\geq & {t-2r+6 \choose r}+(2r-7){t-2r+6 \choose r-1}-(2r-6)\times2^{r-1}({t-2\choose r-2}-1).
\end{eqnarray}
Repeated using the equality ${m+1\choose n}={m\choose n}+{m\choose n-1}$ to the above inequality, we have
$$|[t-1]^{(r)}\bigcap E| \geq {t-1 \choose r}-[(2r-6)\times2^{r-1}+(r-4)(2r-7)]({t-2\choose r-2}-1).$$
So
$$0<|[t-1]^{(r)}\backslash E|\leq [(2r-6)\times2^{r-1}+(r-4)(2r-7)]({t-2\choose r-2}-1).$$
Since $G$ does not contain $[t-1]^{(r)}$.
Let $E^*=E\bigcup [t-1]^{(r)}$ and $G^*=([t], E^*)$. Denote the  number of edges of $G^*$ by $m^*$, then ${t-1\choose r}\leq m^* \leq {t-1\choose r}+ {t-2\choose r-1}-2^{r-3}({t-2\choose r-2}-1).$  So $\lambda(G^*)=\lambda([t-1]^{(r)})$ by Theorem \ref{TheoremPTZ}.
Clearly, $\lambda(G^*,\vec{x})-\lambda(G,\vec{x})>0$  since $x_1\ge x_2 \ge \cdots \ge x_t > 0$ and $|[t-1]^{(r)}\backslash E|>0$. Hence  $\lambda(G)=\lambda(G,\vec{x})<\lambda(G^*,\vec{x})\leq \lambda(G^*)=\lambda([t-1]^{(r)}).$ This completes the proof of Theorem \ref{mainresult}.\qed

%\section{Partial results for  Conjecture \ref{conjecture}  and Conjecture \ref{conjecture2}}\label{section2}

%\begin{coro}\label{coroturan2}  Let $G=(V,E)$ be a left-compressed $3$-graph with  $m$ edges. If $m\geq \frac{(t-3)(t-2)t^3}{6(t-1)^2}$, then $G$ contains a clique order of $\lfloor \frac{t-2r+6}{3}\rfloor$.
%\end{coro}
%\noindent{\bf  Proof} $G=(V,E)$ be a $3$-graph with $t$ vertices and $m$ edges. Assume that $m\geq \frac{(t-3)(t-2)t^3}{6(t-1)^2}$. Clearly, $x_1=x_2=\cdots=x_t=\frac{1}{t}$ is a weighting for $G$. Hence $\lambda(G)\geq \frac{(t-3)(t-2)t^3}{6(t-1)^2}\frac{1}{t^3}=\frac{(t-3)(t-2)}{6(t-1)^2}=\lambda([t-1]^{(r-1)}).$ Hence $G$ contains a clique order of $\lfloor \frac{t-2r+6}{3}\rfloor$ by Theorem \ref{theorem 2}. This completes the proof. \qed

%Let ${\vec y}=(y_1, y_2, \ldots, y_t)$ given by $y_i={1 \over t}$ for each $i, 1\le i\le t$. Then
%$${2 \over 27}\ge \lambda(G)\ge \lambda(G, {\vec y})={m \over t^3}.$$
%So
%$$d(G)={m \over {t\choose 3}}\le {2t^3 \over 27{t\choose 3}}={4t^3 \over 9(t-1)(t-2)(t-3)}.$$ \qed

%It is well-known that the Tur\'an density of $K_4^{(r-1)}$ is at least ${5 \over 9}$. Several decades ago, Tur\'an conjectured that Tur\'an density of $K_4^{(r-1)}$ is exactly ${5 \over 9}$. This conjecture still remains open.  However, for left-compressed $K_4^{(r-1)}$-free 3-uniform graphs, we just showed that  densities are no more than ${4 \over 9}$.

\section{Remarks and conclusions}\label{section3}
We remark that, in the proof of Theorem \ref{Talr}, we see that $\gamma_r=2^{2^r}$ and $t\ge \kappa_0(r)$, where $\kappa_0(r)$ is  a sufficiently large integer such that ${t-2 \choose r-1} \geq\gamma_r (t-1)^{r-2}=2^{2^r}(t-1)^{r-2}$ for $t\ge \kappa_0(r)$. In Corollary \ref{improvement},  we improve the upper bound of $m$ from ${t-1 \choose r} + {t-2 \choose r-1} - \gamma_r (t-1)^{r-2}$ to ${t-1\choose r}+ {t-2\choose r-1}-[(2r-6)\times2^{r-1}+2^{r-3}+(r-4)(2r-7)-1]({t-2\choose r-2}-1).$ Correspondingly, we improve the condition on $t$ from ${t-2 \choose r-1}\geq2^{2^r}(t-1)^{r-2}$ to ${t-2 \choose r-1}\geq[(2r-6)\times2^{r-1}+2^{r-3}+(r-4)(2r-7)-1]({t-2\choose r-2}-1).$

The method developed in the proof of Theorem \ref{mainresult} can also be used to deal with the case for $r=3$ ( see \cite{TPZZ} ). A natural question in the future study is how to prove  similar results as Theorem \ref{mainresult} and Corollary \ref{improvement} without the restriction of  the number of vertices. This will be considered in the future work.

\noindent{\bf Acknowledgments}  This research is partially  supported by National Natural Science Foundation of China (Grant No.61304021
).

\end{document}